\documentclass{amsart}
\usepackage[all,cmtip,ps]{xy}
\usepackage{amsmath}
\usepackage{amsthm}
\usepackage{amssymb}
\usepackage{amscd}
\usepackage[usenames]{color}

\newcommand*{\D}{\mathcal{D}}

\DeclareMathOperator{\Gen}{Gen}

\DeclareMathOperator{\Le}{\mathbb{L}}

\DeclareMathOperator{\Hom}{Hom}
\DeclareMathOperator{\End}{End}
\DeclareMathOperator{\Ext}{Ext}
\DeclareMathOperator{\Tor}{Tor}
\DeclareMathOperator{\Ker}{Ker}

\newcommand*{\K}{\mathcal{K}}

\DeclareMathOperator{\R}{\mathbb{R}}

\def\b#1{{#1}^{\bullet}}

\newcommand*{\rMod}{\textrm{\textup{Mod-}}}

\newtheorem{thm}{Theorem}[section]

\newtheorem{lemma}[thm]{Lemma}
\newtheorem{prop}[thm]{Proposition}
\newtheorem{cor}[thm]{Corollary}
\theoremstyle{definition}

\newtheorem{de}[thm]{Definition}

\def\dualita#1#2{\mathrel{
                 \mathop{\vcenter{
                 \offinterlineskip
                 \hbox to 1.2truecm{\rightarrowfill}
                 \hbox to 1.2truecm{\leftarrowfill}}}%
                 \limits_{#2}^{#1}}}
\def\dual{\mathrel{
                 \mathop{\vcenter{
                 \offinterlineskip
                 \hbox to .6truecm{\rightarrowfill}
                 \hbox to .6truecm{\leftarrowfill}}}%
                 }}

\begin{document}
\title[Derived equivalence]
{Derived equivalence induced by $n$-tilting modules}
\author[S. Bazzoni]{Silvana Bazzoni}
\author[F. Mantese]{Francesca Mantese}
\author[A. Tonolo]{Alberto Tonolo}

\address[Silvana Bazzoni]{Dipartimento di 
Matematica Pura ed Applicata\\ Universit\`a di Padova\\ via Belzoni 7, 
I-35131 Padova - Italy}
\email[Silvana Bazzoni]{bazzoni@math.unipd.it}
\address[F. Mantese]{Dipartimento di Informatica, Universit\`a degli Studi di Verona, strada Le Grazie  15, I-37134 Verona - Italy}
\email{francesca.mantese@univr.it}
\address[Alberto Tonolo]{Dipartimento di 
Matematica Pura ed Applicata\\ Universit\`a di Padova\\ via Belzoni 7, 
I-35131 Padova - Italy}
\email[Alberto Tonolo]{tonolo@math.unipd.it}
\thanks{Research supported by grant CPDA071244/07 of Padova University}
\begin{abstract}
Let $T_R$ be a right $n$-tilting module over an arbitrary associative ring $R$.
In this paper we prove that there exists a
$n$-tilting module $T'_R$ equivalent to $T_R$ which induces a derived equivalence between  the unbounded derived category $\D(R)$ and a triangulated subcategory $\mathcal E_{\perp}$ of $\D(\End(T'))$ equivalent to the quotient category of $\D(\End(T'))$ modulo  the kernel of the total left derived functor $-\otimes^{\mathbb L}_{S'}T'$. In case $T_R$ is a classical $n$-tilting module, we get again the Cline-Parshall-Scott and Happel's results.
\end{abstract}
\maketitle
\section*{Introduction}
Tilting theory generalizes the classical Morita theory of  equivalences between module categories. Originated in the works of Gel'fand and Ponomariev, Brenner and Butler, Happel and Ringel \cite{BeGePo73,BrBu80,HaRi82}, it has been generalized in various directions.
In the recent literature, given an associative ring $R$ with $0\not=1$, a right $R$-module $T_R$ is said to be $n$-tilting if the following conditions are satisfied:
\begin{enumerate}
\item[(T1)] there exists a projective resolution of right $R$-modules
\[0\to P_n\to ...\to P_1\to P_0\to T\to 0;\]
\item[(T2)] $\Ext^i_R(T, T^{(\alpha)})=0$ for each $i>0$ and each cardinal $\alpha$;
\item[(T3)] there exists a coresolution of right $R$-modules
\[0\to R\to  T_0\to T_1\to ...\to T_m\to 0,\]
where the $T_i$'s are direct summands of arbitrary direct sums of copies of $T$.
\end{enumerate}
If the projectives $P_i$'s in (T1) can be assumed finitely generated, then the $n$-tilting module $T_R$ is said \emph{classical $n$-tilting}. 

Let us denote by $S=\End(T_R)$ the endomorphism ring of $T$ and by $KE_i(T)$ and $KT_i(T)$, $0\leq i\leq n$, the following classes
\[KE_i(T)=\{M\in\rMod R:\Ext^j_R(T,M)=0\text{ for each }0\leq j\not=i\},\]\[
KT_i(T)=\{N\in\rMod S:\Tor_j^S(N,T)=0\text{ for each }0\leq j\not=i\}.\]

In 1986 Miyashita \cite{Mi86} proved that if $T_R$ is a classical $n$-tilting, then the functors
$\Ext^i_R(T,-)$ and $\Tor_i^S(-,T)$ induce equivalences between the classes $KE_i(T)$ and $KT_i(T)$.

In the same years, works of several authors showed that the natural context for studying equivalences induced by classical tilting modules is that of derived categories. In particular Cline, Parshall and Scott \cite{ClPaSc86}, generalizing a result of Happel \cite{Ha87}, proved that a classical $n$-tilting module $T_R$ provides a derived equivalence between the bounded derived categories $\D^b(R)$ and $\D^b(S)$ of bounded cochain complexes of right $R$- and $S$- modules.

In 1988 Facchini \cite{Fa87,Fa88} proved that, over a commutative domain, the divisible module $\partial$ introduced by Fuchs \cite{Fu84} is an infinitely generated 1-tilting module and it provides a pair of equivalences 
\[KE_0(\partial)\dualita{\Hom(\partial,-)}{-\otimes\partial}KT_0(\partial)\cap I\text{-}Cot,\quad
KE_1(\partial)\dualita{\Ext^1(\partial,-)}{Tor_1(-,\partial}KT_1(\partial)\cap I\text{-}Cot
\]
between the category $KE_0(\partial)$ of all divisible modules and the category $KT_0(\partial)\cap I\text{-}Cot$ of all \emph{$I$-reduced $I$-cotorsion modules}, and the category $KE_1(\partial)$ of all reduced modules and the category $KT_1(\partial)\cap I\text{-}Cot$ of all \emph{$I$-divisible $I$-cotorsion} modules, respectively.
In 1995 Colpi and Trlifaj \cite{CoTr95} started the study in general of $1$-tilting modules. They realized that it can be useful to ``change slightly'' the tilting module to realize a good equivalence theory.   They proved that if $T_R$ is a 1-tilting module, there exists another 1-tilting module $T'_R$ \emph{equivalent to} $T_R$ (i.e. $KE_0(T)=KE_0(T')$), with endomorphism ring $S'=\End(T')$, such that the functors $\Hom_R(T',-)$ and $-\otimes_{S'}T'$ induce an equivalence between $KE_0(T)=KE_0(T')$ and its image class in $\rMod S'$. Moreover $T'$ results to be a finitely presented $S'$-module.
In 2001 Gregorio and Tonolo extended this result proving the existence of a pair of equivalences
\[KE_i(T')\dualita{\Ext^i_R(T',-)}{\Tor_i^{S'}(-,T')}KT_i(T')\cap Cost(T'),\quad i=1,2\]
where $Cost(T')$ is the class of \emph{costatic} right $S'$-modules (see \cite{GrTo01}).

In 2009 Bazzoni \cite{Ba09} gives a better understanding of the whole situation in the setting of derived categories proving that for a $1$-tilting module $T_R$ it is possible to find an equivalent $1$-tilting module $T'$ which induces a derived equivalence between  the unbounded derived category $\D(R)$ and the quotient category of $\D(S')$ modulo the full triangulated subcategory $\Ker(-\otimes_{S'}^{\mathbb L}T')$, namely the kernel of the total left derived functor of the functor $-\otimes_{S'}T'$.

In this paper we generalize the Bazzoni's result to a general $n$-tilting module $T_R$. We prove the existence of a \emph{good} $n$-tilting module $T_R'$ \emph{equivalent} to $T_R$ (see Definition~\ref{def:equivalent}) which, also in such a case, provides a derived equivalence between  the unbounded derived category $\D(R)$ and a triangulated subcategory $\mathcal E_{\perp}$ of $\D(\End(T'))$. The category $\mathcal E_{\perp}$ results to be equivalent to the quotient category of $\D(\End(T'))$ modulo  the kernel of the total left derived functor $-\otimes^{\mathbb L}_{S'}T'$. Moreover, as done in \cite{MaTo09} in the contravariant case, we interpret the derived equivalence at the level of stalk complexes obtaining on the underlying module categories a generalization of the Miyashita equivalences.

\section{$n$-tilting classes}

In 2004 Bazzoni (see \cite{Ba04}) proved that $T_R$ is a $n$-tilting module if and only if the classes
\[T^{\perp_\infty}:=\{M_R:\Ext^i_R(T,M)=0 \text{ for each }i>0\}\]
and
\[\Gen_n(T):=\{M_R:\ \exists\ 
T^{(\alpha_n)}\to ...\to T^{(\alpha_1)}\to M\to 0, \text{ for some cardinals }\alpha_i\}
\]
coincide.
\begin{de}\label{def:equivalent}
Two $n$-tilting right $R$-modules $T_R$ and $T'_R$ are said \emph{equivalent} if $\Gen_n(T_R)=\Gen_n(T'_R)$.
\end{de}
An arbitrary direct sum of copies of a $n$-tilting module is a $n$-tilting module equivalent to the original one. Therefore
equivalent tilting modules can have completely different endomorphism rings. 
\begin{de}\label{def:goodtilting}
We say that $T_R$ is a \emph{good} $n$-tilting module if it is $n$-tilting and it satisfies the condition
\begin{enumerate}
\item[(T3')] there is an exact sequence
\[0\to R\to T_0\to T_1\to ...\to T_n\to 0\]
where the $T_i$'s are direct summands of finite direct sums of copies of $T$.
\end{enumerate}
\end{de}
Each classical $n$-tilting module is good \cite[Section 5.1]{GoTr06}.
\begin{prop}
For any $n$-tilting module $T_R$ there exists an equivalent good $n$-tilting module $T'_R$ such that
\[KE_i(T)=KE_i(T') \text{ for each }i\geq 0
.\]
\end{prop}
\begin{proof}
Let $T_R$ be a $n$-tilting module.
If it is classical, then $T$ already satisfies (T3'). Otherwise,
from condition $(T3)$
we easily get the exact sequence
\[0\to R\to T_0\to ...\to T_{n-2}\to T_{n-1}\oplus T_n^{(\omega)}\to T_n\oplus T_n^{(\omega)}\to 0\]
that can be rewritten in the form
\[0\to R\to T_0\to ...\to T_{n-2}\to T_{n-1}\oplus T_n^{(\omega)}\to T_n^{(\omega)}\to 0.\]
With the same argument we get the exact sequence
\[0\to R\to  ...\to T_{n-3}\to T_{n-2}\oplus (T_{n-1}\oplus T_n^{(\omega)})^{(\omega)}\to T_{n-1}\oplus T_n^{(\omega)}\oplus
(T_{n-1}\oplus T_n^{(\omega)})^{(\omega)}\to T_n^{(\omega)}\to 0,\]
and hence the exact sequence
\[0\to R\to T_0\to ...\to T_{n-3}\to T_{n-2}\oplus T_{n-1}^{(\omega)}\oplus T_n^{(\omega)}\to T_{n-1}^{(\omega)}\oplus T_n^{(\omega)}\to
T_n^{(\omega)}\to 0.\]
Iterating this procedure we get an exact sequence
\[0\to R\to T_0\oplus T_1^{(\omega)}\oplus ...\oplus T_n^{(\omega)}\to ...\to T_{n-2}^{(\omega)}\oplus T_{n-1}^{(\omega)}\oplus T_n^{(\omega)}\to T_{n-1}^{(\omega)}\oplus T_n^{(\omega)}\to
T_n^{(\omega)}\to 0.\]
Let $T'=T_0\oplus T_1^{(\omega)}\oplus ...\oplus T_n^{(\omega)}$;
since $T'$ is a direct summand of a direct sum of copies of $T$, we have
\[\Gen_n(T')\subseteq \Gen_n(T)=T^{\perp_{\infty}}\subseteq T'^{\perp_{\infty}},\]
and $T'$ satisfies properties $(T1)$ and $(T2)$ of tilting modules. Since by construction it satisfies also property (T3'), we have $\Gen_n(T')=T'^{\perp_{\infty}}$ and $T'$ is the wanted good $n$-tilting equivalent to $T$. 

Finally, since $\Ker\Ext^j(T,-)=\Ker\Ext^j(T_0\oplus ... \oplus T_n,-)=\Ker\Ext^j(T',-)$, we conclude that $KE_i(T)=KE_i(T')$ for each $i\geq 0$.
\end{proof}
A good $n$-tilting module has an endomorphism ring $S$ sufficiently large to permit to build a good equivalence theory between the unbounded derived categories $\D(R)$ and $\D(S)$. In the sequel we will work directly with good $n$-tilting modules.


\begin{prop}
Let $T_R$ be a good $n$-tilting module and $S=\End(T_R)$. Then ${}_ST$ has a projective resolution
\[0\to Q_n\to ...\to Q_0\to {}_ST\to 0\]
where the $Q_i$'s are direct summand of a finite direct sum of copies of $S$, $\Ext^i_S(T,T)=0$ for each $i\geq 0$, and $R\cong \End({}_ST)$.
\end{prop}
\begin{proof}
By Definition~\ref{def:goodtilting} there is an exact sequence
\[0\to R\to T_0\to T_1\to ...\to T_n\to 0\]
with the $T_i$'s direct summands of $T^m$ for a suitable $m\in\mathbb N$. Denote by $K_i$ the kernel of the map $T_i\to T_{i+1}$, $1\leq i\leq n-1$. Applying the contravariant functor $\Hom_R(-, T)$ we get easily by dimension shifting that
\[0=\Ext^i_R(K_j,T) \text{ for each }1\leq j\leq n-1,\text{ and }i\geq 1.\]
Therefore we have the exact sequence
\[(\dag)\quad 0\to \Hom_R(T_n,T)\to \Hom_R(T_{n-1},T)\to ...\to \Hom_R(T_{1},T)\to \Hom_R(T_{0},T)\to {}_ST\to 0\]
where each $\Hom_R(T_i,T)$ is a direct summand of $\Hom_R(T^m,T)=S^m$ and hence a finitely generated projective $S$-module. Given a right $R$-module $M$, let us denote for semplicity by $M^*$ the left $S$-module
$\Hom_R(M,T)$, by $M^{**}$ the right $R$-module $\Hom_S(M^*, T)$,  and by $\delta_M$ the evaluation map $M\to M^{**}$. The modules $K^*_i$ are the cokernels of the morphisms $\Hom_R(T_{i+1},T)\to \Hom_R(T_{i},T)$, $1\leq i\leq n-1$. Applying to $(\dag)$ the contravariant functor $\Hom_S(-, T)$ we get the following commutative diagrams with exact rows:
\[
\xymatrix{
0\ar[r]&\Hom_S(T,T)=R^{**}\ar[r]&T_0^{**}\ar[r]&K_1^{**}\ar[r]&\Ext^1_S(T,T)\ar[r]&0\\
0\ar[r]&R\ar[u]^{\delta_R}\ar[r]&T_0\ar[u]^{\delta_{T_0}}\ar[r]&K_1\ar[u]^{\delta_{K_1}}\ar[r]&0
}
\]
\[\hdots\]
\[\xymatrix{
0\ar[r]&K_{n-1}^{**}\ar[r]&T_{n-1}^{**}\ar[r]&T_n^{**}\ar[r]&\Ext^1_S(K_{n-1}^*,T)\ar[r]&0\\
0\ar[r]&K_{n-1}\ar[u]^{\delta_{K_{n-1}}}\ar[r]&T_{n-1}\ar[u]^{\delta_{T_{n-1}}}\ar[r]&T_n\ar[u]^{\delta_{T_n}}\ar[r]&0
}
\]
Since the $\delta_{T_i}$'s are isomorphisms we get
\[\Ext^1_S(T,T)=0\text{ and } 0=\Ext^1_S(K^*_{i},T)\cong\Ext^{i+1}_S(T,T) \text{ for each }1\leq i\leq n-1,\]
and $R\cong\Hom_S(T,T)$.

\end{proof}

\begin{lemma}[Lemmas~1.8, 1.9 \cite{Mi86}]\label{lemma:Miya}
Let $T_R$ be a good $n$-tilting and $S=\End T$. For any right $R$-module $M$ in $T^{\perp_\infty}$ and any right projective $S$-module $P_S$, we have
\begin{enumerate}
\item $\Tor^S_i(\Hom_R(T,M),T)=0$ for each $i>0$.
\item $\Hom_R(T,M)\otimes_S T\cong M,\quad f\otimes t\mapsto f(t)$
\item $\Ext^i_R(T,P\otimes_S T)=0$ for each $i>0$.
\end{enumerate}
If $T_R$ is a classical $n$-tilting module, then
\begin{enumerate}
\item[(4)] $P\cong \Hom_R(T,P\otimes_S T),\quad p\mapsto (f:t\mapsto p\otimes t)$.
\end{enumerate}
\end{lemma}
\begin{proof}
Everything except condition (3) follows by the quoted lemmas in \cite{Mi86}. If $P\leq^{\oplus} S^{(\alpha)}$ we have
\[\Ext^i_R(T,P\otimes_S T)\leq^{\oplus}\Ext^i_R(T,S^{(\alpha)}\otimes_S T)=\Ext^i_R(T,T^{(\alpha)})=0.\]
\end{proof}

\section{Tilting equivalences in derived categories}
In the sequel, for any ring $R$, we denote by $\K(R)$   the homotopy category of unbounded complexes of right $R$-modules  and by $\D(R)$  the associated derived category. 
Given an object $M\in\rMod R$, we  continue to denote by $M$ also the \emph{stalk complex} in $\D(R)$ associated to $M$, i.e. the complex with $M$ concentrated in degree zero.
Any complex $C^{\bullet}\in\D(R)$ admits a  $K$-injective resolution, i.e. a complex $\underline{\textbf i}C^{\bullet}$ quasi-isomorphic to $C^{\bullet}$ whose terms are injective modules. Similarly, any complex $C^{\bullet}\in\D(R)$ admits a  $K$-projective resolution, i.e. a complex $\underline{\textbf p}C^{\bullet}$ quasi-isomorphic to $C^{\bullet}$ whose terms are projective modules (see for instance \cite{BoNe93}).  This result guarantees the existence of the total derived functor of any additive functor defined on module categories.

Given any covariant left exact functor $H : \rMod R\to \rMod S$, we denote by $\R H$ its total right derived functor defined on $\D(R)$. For any $C^{\bullet}\in \D(R)$, $\R H(C^{\bullet})$ coincides with the complex $H(\underline{\textbf i}C^{\bullet})$, where we still denote by $H$ its extension to $\K(R)$. Similarly, for any right exact covariant functor $G: \rMod S\to \rMod R$, we denote by $\Le G$ its total left derived functor defined on $\D(S)$. For any $N^{\bullet}\in \D(S)$, $\Le G(N^{\bullet})$ coincides with the complex $G(\underline{\textbf p}N^{\bullet})$. 

A module  $M$ in $\rMod R$ is called \emph{$H$-acyclic} if $R^i H M:=H^i(\R H M)=0$ for any $i\not=0$. The abelian group  $R^iH M$ coincides with the usual $i$-th derived functor $H^{(i)}(-)$ of $H$ evaluated in $M$. Analogously \emph{$G$-acyclic} objects are defined and
$L^i G(-):=H^i(\Le G (-))=G^{(-i)}(-)$. In view of these consideration, by Lemma~\ref{lemma:Miya} we have immediately the following result.

\begin{cor}\label{cor:aciclici}
Let $T_R$ be a good $n$-tilting module with endomorphism ring $S$. Then for each injective module $I_R$ and each projective module $P_S$ we have
\begin{enumerate}
\item $\Hom_R(T,I)$ is $-\otimes_S T$-acyclic;
\item $P\otimes_S T$ is $\Hom_R(T,-)$-acyclic.
\end{enumerate}
In particular for cochain complexes $\b I$ and $\b P$ whose terms are injective right $R$-modules and projective right $S$-modules respectively, we have
\[\R\Hom(T, \b I)\otimes_S^{\Le}T=\Hom(T, \b I)\otimes_ST\text{ and }
\R\Hom(T, \b P\otimes_S^{\Le}T)=\Hom(T, \b P\otimes_ST).
\]
\end{cor}
%


Finally,  we recall that any adjoint pair of functors $(G, H)$ between categories of modules induces an adjoint pair $(\Le G, \R H)$ between the associated unbounded derived categories.
For other notations and results in derived categories we refer to \cite{Ha66,We94}.

In the sequel we denote by $H$ the functor $\Hom_R(T,-)$ and by $G$ the functor $-\otimes_S T$.

\begin{thm}
Let $T_R$ be a good $n$-tilting module and $S=\End T_R$. The following hold:
\begin{enumerate}
\item The counit adjunction morphism
\[\mathbb LG\circ \mathbb R H\to Id_{\D(R)}\]
is invertible;
\item the functor $\mathbb R H:\D(R)\to\D(S)$ is fully faithful;
\item if $\Sigma$ is the system of morphisms $u\in\D(S)$ such that $\mathbb LGu$ is invertible in $\D(R)$, then $\Sigma$ admits a calculus of left fractions and the category $\D(S)[\Sigma^{-1}]$ coincides with the quotient category $\D(S)$ modulo the full triangulated subcategory $\Ker(\mathbb LG)$ of the objects annihilated by the functor $\mathbb L G$;
\item there is a triangle equivalence
\[\D(S)[\Sigma^{-1}]\dualita{\Theta}{\mathbb R H} \D(R)\] where $\Theta$ is the functor such that $\mathbb LG=\Theta\circ q$ with $q$ the canonical quotient functor $q:\D(S)\to \D(S)[\Sigma^{-1}]$.
\end{enumerate}
\end{thm}
\begin{proof}
(1) Let $M^{\bullet}$ be a complex in $\D(R)$ and consider a $K$-injective resolution $\underline{\textbf i}M^{\bullet}$ of $M^{\bullet}$. 
By Corollary~\ref{cor:aciclici} we have
\[\mathbb LG(\mathbb RH(M^{\bullet}))=\mathbb LG(H(\underline{\textbf i}M^{\bullet}))=G(H(\underline{\textbf i}M^{\bullet})).\]
Since $(\Hom_R(T,I^n)\otimes_S T)_{n\in\mathbb Z}$ and $\underline{\textbf i}M^{\bullet}$ are isomorphic by Lemma~\ref{lemma:Miya}, (2), we have
\[\mathbb LG(\mathbb RH(M^{\bullet}))=G(H(\underline{\textbf i}M^{\bullet}))\cong \underline{\textbf i}M^{\bullet}= M^{\bullet}.\]
Conditions (2), (3) and (4) follow by applying \cite[Proposition~I.1.3]{GaZi67}.
\end{proof}


Let $\mathcal C$ be a triangulated category closed under arbitrary coproducts; recall that a triangle functor $L:\mathcal C\to \mathcal C$ is a \emph{Bousfield localization}  if there exists a natural transformation $\phi:1_{\mathcal C}\to L$ such that for each $X$ in $\mathcal C$
\begin{enumerate}
\item $L(\phi_X):L(X)\to L^2(X)$ is an isomorphism;
\item $L(\phi_X)=\phi_{L(X)}$.
\end{enumerate}
In such a case the kernel $\mathcal L$ of $L$ is a full triangulated subcategory of $\mathcal C$ closed under coproducts, i.e. it is a \emph{localizing} subcategory. The category 
\[\mathcal L_{\perp}:=\{X\in\mathcal C: \Hom_{\mathcal C}(\mathcal L,X)=0\}\]
is called the subcategory of $\mathcal L$-local objects. If also $\mathcal L_{\perp}$ is closed under coproducts, then $\mathcal L$ is called \emph{smashing} \cite{Bo79,BoNe93}.

A localization functor ${L}$ factorizes as
\[\mathcal C\stackrel{q}{\to}\mathcal C/\Ker L\stackrel{\rho}{\underset{\cong}{\longrightarrow}} \mathcal L_{\perp}\overset{j}{\hookrightarrow}\mathcal C\]
where $q$ is the canonical quotient functor and $\rho$ is an equivalence; $(\rho\circ q, j)$ is an adjoint pair. Moreover the composition
\[ \mathcal L_{\perp}\overset{j}{\hookrightarrow}\mathcal C
\stackrel{q}{\to}\mathcal C/\Ker L\]
is an equivalence and $(q,j\circ\rho)$ is an adjoint pair (see \cite[Section~4]{BoNe93}, or \cite[Proposition~1.6]{AlLoSo00}, or \cite[Propositions 4.9.1, 4.11.1]{Kr09}).

\begin{thm}\label{teo:generale}
Let $(\Phi, \Psi)$ be an adjoint pair of covariant functors between 
triangulated categories
\[\mathcal C\dualita{\Phi}{\Psi}\mathcal D.\]
Denote by 
$\phi:1_{\mathcal C}\to \Psi\circ\Phi$ and 
$\psi:\Phi\circ\Psi\to1_{\mathcal D}$ the corresponding unit and counit.
If $\psi$ is a natural isomorphism, then the functor $L:=\Psi\circ\Phi$ is a localization functor with kernel $\mathcal L=\Ker\Phi$. The functor $\Psi$ factorizes through $\mathcal L_{\perp}$ as $\Psi=j\circ\overline{\Psi}$, where $j$ is the inclusion $\mathcal L_{\perp}{\hookrightarrow}\mathcal C$. Finally we have a triangle equivalence
\[
\mathcal L_{\perp}\dualita{\Phi\circ j}{\overline{\Psi}}\mathcal D
\]
where $\Phi\circ j$ is the restriction of $\Phi$ to $\mathcal L_{\perp}$ and $\overline{\Psi}$ is the corestriction of $\Psi$ to $\mathcal L_{\perp}$.
\end{thm}
\begin{proof}
Since $(\Phi, \Psi)$ is an adjoint pair, we have 
\[\psi_{\Phi(X)}\circ\Phi(\phi_X)=1_{\Phi(X)};\]
applying the functor $\Psi$ we get
\[\Psi(\psi_{\Phi(X)})\circ L(\phi_X)=1_{L(X)}.\]
On the other hand, again by the adjunction, we have
\[\Psi(\psi_{\Phi(X)})\circ\phi_{\Psi\Phi(X)}=1_{\Psi\Phi(X)},\text{ i.e. }\Psi(\psi_{\Phi(X)})\circ\phi_{L(X)}=1_{L(X)}.\]
Since $\psi_{\Phi(X)}$ is an isomorphism by assumption, we have that for each $X$ in $\mathcal C$
\[L(\phi_X)=\phi_{L(X)}=(\Psi(\psi_{\Phi(X)}))^{-1}\]
is an isomorphism. Hence $L$ is a localization functor.

An object $X$ belongs to $\mathcal L=\Ker L$ if and only if we have $0=\Phi(0)=\Phi(\Psi\Phi(X))\cong \Phi(X)$.

Next, since $L=\Psi\circ\Phi$ factorizes through $\mathcal L_{\perp}$ and $\Phi(\Psi(Y))\cong Y$ for each $Y$ in $\mathcal D$, also $\Psi$ factorizes through $\mathcal L_{\perp}$. Therefore we have the following commutative diagram:
\[
\xymatrix{
{}\mathcal L_{\perp}\ar@{^(->}
[r]|j\ar@/^2pc/[rr]^{q\circ j}_{\cong}&\mathcal C\ar[r]|-q\ar[drr]|{\Phi}&\mathcal C/\Ker \Phi\ar[rr]|{\rho}_{\cong}\ar[dr]|{\Theta}&&\mathcal L_{\perp}\ar@{^(->}
[r]|j&\mathcal C\\
&&&\mathcal D\ar[ur]|{\overline{\Psi}}\ar[urr]|{\Psi}
}
\]
Finally $\Phi\circ j\circ\overline{\Psi}=\Phi\circ\Psi\cong 1_{\mathcal D}$, and
$\overline{\Psi}\circ\Phi\circ j=\rho\circ q\circ j$, being a composition of two equivalences, is naturally isomorphic to $1_{\mathcal L_{\perp}}$.
\end{proof}

Applying Theorem~\ref{teo:generale} to our context we obtain the following result
\begin{cor}\label{teo:trequi}
Let $T_R$ be a good $n$-tilting $R$-module and $S=\End(T)$. Denoted by $\mathcal E$ the kernel of $\Le G$, and denoting by $\R H$ and $\Le G$ also their restriction and corestriction, we have a triangulated equivalence
\[\D(R)\dualita{\mathbb R H}{\mathbb L G}\mathcal E_{\perp}.\]
\end{cor}

Embedding right $R$-modules and $S$-modules in $\D(R)$ and $\D(S)$ via the canonical functor, we obtain the following generalization of the Miyashita's results \cite[Theorem~1.16]{Mi86}:

\begin{cor}
Let $T_R$ be a good $n$-tilting $R$-module and $S=\End(T)$.
Then for each $0\leq i\leq n$ there is an equivalence
\[KE_i\dualita{\Ext^i_R(T,-)}{\Tor^S_i(-,T)} KT_i\cap \mathcal E_\perp\]
\end{cor}
\begin{proof}
Let $M\in KE_i$; then by Corollary~\ref{teo:trequi}, $\R H(M)=R^i H(M)[-i]=\Ext^i_R(T,M)[-i]$ belongs to $\mathcal E_{\perp}$. Since $\mathcal E_{\perp}$ is closed under shift, $\Ext^i_R(T,M)\in \mathcal E_{\perp}$.
In $\D(R)$, by Theorem~\ref{teo:trequi}, (1), we have
\[M\cong \Le G\R H(M)=\Le G(\Ext^i_R(T,M)[-i]);
\]
then for each $j\not=0$
\[0=H^j\Le G(\Ext^i_R(T,M)[-i])=H^{j-i}\Le G(\Ext^i_R(T,M))=\Tor^S_{i-j}(\Ext^i_R(T,M),T).
\]
Therefore $\Ext^i_R(T,M)$ belongs to $KT_i\cap \mathcal E_\perp$ and $M\cong \Tor^S_{i}(\Ext^i_R(T,M),T)$. 
Analogously if $N\in KT_i\cap \mathcal E_\perp$, then 
\[\Le G(N)=L^{-i} G(N)[i]=\Tor^S_i(N,T)[i]\]
and since $\R H\Le G(N)=N$ in $\D(S)$, necessarily $\Tor^S_i(N,T)$ belongs to $KE_i$ and
$N\cong \Ext^i_R(T,\Tor^S_{i}(N,T))$. 
\end{proof}


\begin{prop}
The following are equivalent:
\begin{enumerate}
\item $T_R$ is a classical $n$-tilting;
\item $\mathcal E=0$ or equivalently $\mathcal E_{\perp}=\D(S)$;
\item the class $\mathcal E$ is smashing.
\end{enumerate}
\end{prop}
\begin{proof}
($1\Rightarrow 2$). Let $\b N$ be a complex in $\mathcal E$ and $\underline{\textbf{p}}\b N$ a $K$-projective resolution of $\b N$. By Lemma~\ref{lemma:Miya}, (3) and (4), we have
\[0=\mathbb RH(\mathbb L G\b N)=\mathbb RH(\mathbb L G\underline{\textbf{p}}\b N)=\mathbb RH(\underline{\textbf{p}}\b N\otimes_S T)=\]
\[=\Hom_R(T,\underline{\textbf{p}}\b N\otimes_S T)\cong \underline{\textbf{p}}\b N=\b N.\]
We conclude that $\mathcal E=0$ by Corollary~\ref{teo:trequi}. 

($2\Rightarrow 3$) is obvious.

($3\Rightarrow 2$). 
Since $S=\R H(T_R)$,  $\mathcal E_{\perp}$ contains the bounded complexes of finitely generated projective $S$-modules,  that is $\mathcal E_{\perp}$ containe the set $\mathcal T^c$ of the compact objects of $\D(S)$.

Since $\D(S)$ is compactly generated by $\mathcal T^c$, $\D(S)$ is the smallest triangulated category closed under coproducts and containig $\mathcal T^c$.
Thus, if $\mathcal E_{\perp}$ is closed under coproducts we get that $\mathcal E_{\perp}=\D(S)$, hence $\mathcal E=0$

($2\Rightarrow 1$). By \cite[Propositions ~6.2, 6.3 and Theorem 6.4]{Ri89} for any equivalence \[\D^b(R)\dualita{\Psi}{\Phi}\D^b(S)\] it is $\Psi=\R \Hom(\Phi(S), -)$ and $\Phi=-\otimes^{\mathbb L}_S\Psi(R)$ with $\Phi(S)$ isomorphic to a bounded complex of finitely generated projective $R$-modules. Since 
\[\Le G(S)=G(S)=S\otimes T=T_R,\]
we conclude that $T_R$ is a classical $n$-tilting module.
\end{proof}


\end{document}